\documentclass {article}
\title{Theta-functions on noncommutative tori.} 
\author {Albert Schwarz \thanks{Partially  supported by  NSF grant 
No. DMS 9801009}} 
\begin {document}
\maketitle
\begin {abstract}
Ordinary theta-functions can be considered as holomorphic sections of line
bundles over tori. We  show that one can define generalized
theta-functions as holomorphic elements of projective modules over
noncommutative tori (theta-vectors). The theory of these new objects is
not only more general, but also much simpler than the theory of ordinary
theta-functions. It seems that the theory of theta-vectors should be
closely related to Manin's theory of quantized theta-functions, but we
don't analyze this relation.
\end {abstract}
Theta-function of $g$ variables $\vec {z}=(z^1,...,z^g)$ is defined by the
formula 
\begin {equation}
\vartheta (\vec {z},\Omega)=\sum _{\vec {n}\in {\bf Z}^g} \exp (\pi i \vec
{n}^t\Omega \vec {n} +
2\pi i \vec {n}^t\vec {z})
\end {equation}
where $\Omega$ belongs to Siegel upper-half-space ${\cal H}_g$
(i.e. $\Omega$ is a symmetric  $g\times g$  complex matrix whose imaginary
part is positive definite). It is a holomorphic function on ${\bf
C}^g\times {\cal H}_g$ and obeys 
\begin {equation}
\vartheta (\vec {z}+\vec {m},\Omega)=\vartheta (\vec {z},\Omega),
\end {equation}
\begin {equation}
\vartheta (\vec {z}+\Omega\vec {m},\Omega)=\exp (-\pi i \vec
{m}^t\Omega\vec {m}-2\pi i \vec {m}^t\vec {z}) \vartheta (\vec {z},\Omega)
\end {equation}
for all $\vec {m}\in {\bf Z}^g$.

It follows from (2), (3) that $\vartheta (\vec {z},\Omega)$ can be
considered as a holomorphic section of a line bundle over torus ${\bf
C}^g/{\bf Z}^{2g}$. The function $\vartheta (\vec {z},\Omega)$ can be
considered also as a modular form with respect to the action of ${\rm Sp}
(2g,{\bf Z})$ on both the variables $\vec {z},\Omega$. More precisely, 
$$\vartheta (((C\Omega+D)^t)^{-1}\cdot  \vec
{z},(A\Omega+B)  (C\Omega+D)^{-1})=\zeta _{\gamma}\cdot \det (C
\Omega+D)^{1/2}\exp [\pi i \vec {z}^t \cdot (C\Omega +D)^{-1}C\cdot \vec
{z})] \vartheta (\vec {z},\Omega),$$
where $\zeta _{\gamma}$ is an $8$-th root of $1$ and 
\begin {equation}
\gamma= \left ( \begin {array} {cc}
                                                            A & B \\
                                                            C & D
 \end {array}        \right ) \in \Gamma _{1,2}.
\end {equation}
Here $\Gamma_{1,2}$ denotes the subgroup of ${\rm Sp}(2g,{\bf  Z})$
generated by elements of the form 
\begin {equation}
\gamma= \left ( \begin {array} {cc}
                                                            A & 0 \\
                                                            0 & (A^{-1})^t
 \end {array}        \right ) , \ \  
\gamma= \left ( \begin {array} {cc}
                                                            1 & B \\
                                                            0 & 1
 \end {array}        \right ) , \ \  
\gamma= \left ( \begin {array} {cc}
                                                            0 & -1 \\
                                                            1 & 0
 \end {array}        \right ) ,
\end {equation}
where $ A\in {\rm GL}(g,{\bf Z}),\ \ \  B$ is symmetric with even diagonal
entries.

Our goal is to develop the theory of theta-functions over  noncommutative
tori.   Noncommutative tori are the simplest and most important
noncommutative spaces; the differential geometry of these spaces was
studied thoroughly in many papers starting with the seminal paper [1]. We
will use freely  the constructions and results of papers [1]-[8].  By
definition an algebra $T_{\theta}^d$ of smooth functions on noncommutative
torus is an associative algebra of smooth functions on ordinary  torus
$T^d={\bf  R}^d /{\bf Z}^d$ equipped with star-product 
\begin {equation}
(f\star g)(\vec {x})=\int f(\vec {x}+\theta \vec {u})g(\vec {x} +\vec
{v}) \exp (2\pi i \vec {u} \vec {v})d \vec {u} d  \vec {v}, 
\end {equation}
where $\vec {v}\in {\bf R}^d,\ \ \vec {u}\in ({\bf R}^d)^{\ast}$, and
$\theta :({\bf R}^d)^{\ast} \rightarrow  {\bf R}^d$ is a linear operator
determined by an antisymmetric matrix $\theta ^{\alpha \beta}$. (From now
on we will identify $({\bf R}^d)^{\ast}$ with ${\bf R}^d$ using inner
product $\vec {u}\cdot \vec {v}=\vec {u}^t \vec {v}$.) Alternatively using
Fourier representation we consider an element of $T_{\theta}^d$ as a
series 
$$f=\sum _{\vec {n}\in {\bf Z}^d}f_{\vec {n}}U_{\vec {n}} $$  where
$f_{\vec {n}}$ tends to zero faster than any power and the multiplication
is specified by the formula $$U_{\vec {n}}U_{\vec {m}}=e^{\pi i\vec
{n}\theta \vec {m}}U_{\vec {n}+\vec {m}}.$$ 

Elements $U_{\vec {n}}$ can be considered as additive generators of
$T_{\theta}^d$. One can obtain multiplicative generators of $T_{\theta}^d$
considering only $U_{\alpha}=U_{\vec {e}_{\alpha}},\ \  \alpha =1,...,d$,
where $\vec {e}_1,...,\vec {e}_d$ stands for the standard basis of ${\bf
R}^d$. The elements $U_1,...,U_d$ obey 
\begin {equation}
U_{\alpha}\cdot U_{\beta}=e^{2\pi
i\theta^{\alpha\beta}}U_{\beta}U_{\alpha}.
\end {equation}
This means that we can construct a unitary module over $T_{\theta}^d$
specifying unitary operators $U_1,...,U_d$ obeying (7). (We consider
complex conjugation of functions as an involution on the algebra
$T_{\theta}^d$, elements $U_n$ are unitary with respect to this
involution. All $T_{\theta}^d$-modules we consider are unitary
modules.) Space of continuous (or smooth) sections of a vector bundle over
compact manifold $M$ can be considered as a projective module over
commutative algebra of continuous (respectively, smooth) functions on
$M$. Therefore in noncommutative geometry  projective modules over
associative algebras are considered as analogs of vector bundles.   

Projective modules over  noncommutative torus can be constructed in the
following simple way. Let us consider an irreducible representation of
canonical commutation relations 
\begin {equation}
[\nabla _{\alpha},\nabla _{\beta}]=2\pi i\omega _{\alpha \beta},\ \ \
\nabla _{\alpha}^{\star}=-\nabla _{\alpha}. 
\end {equation}
Here $\omega _{\alpha \beta}$ stands for nondegenerate antisymmetric
matrix. Representing it in block diagonal form with $g$ two-dimensional
blocks ($d=2g$) we see that the relations (8) can be rewritten as 
$$[\nabla _{\alpha}, \nabla _{\alpha +g}]=-[\nabla _{\alpha +g}, \nabla
_{\alpha}]=2\pi i\ \ {\rm  if}\ \    0\leq \alpha <g;$$
 \begin {equation}
 [\nabla _{\alpha},\nabla _{\beta}]=0\ \  {\rm  in\ \  all\ \   other\ \
cases.}
\end {equation} 
More precisely, using linear change of variables from $z^1,...,z^d$ to
$x^1,...,x^d$ we can treat $\omega=\omega _{\alpha
\beta}dz^{\alpha}dz^{\beta}$ as the standard symplectic form.
 Then the representation at hand has a familiar realization by means of
operators 
\begin {equation}
\nabla _{\alpha}=2\pi i\hat {x}^{\alpha}\ \ \ {\rm for}\   1\leq \alpha
\leq g,\ \ \ \ \nabla _{\alpha}=\partial _{\alpha}\ \ \  {\rm for}\
g<\alpha \leq 2g, 
\end {equation}
where operators $\hat {x}^{\alpha}$ and $\partial _{\alpha}$ act on the
space ${\cal S}({\bf R}^g)={\cal S}$ of Schwartz functions depending on
$(x^1,...,x^g)$ as operators of multiplication by $x^{\alpha}$ and
differentiation with respect to $x^{\alpha}$.

It is easy to check that the operators 
\begin {equation}
 U_{\alpha}=e^{-\theta^{\alpha \beta} \nabla _{\beta}}
\end {equation}
where $\theta $ stands for the matrix $\theta ={\omega}^{-1}$, specify a
projective module over $T_{\theta}^d$. 

The operators  $\nabla _{\alpha}$ obey the relation  
$$[\nabla _{\alpha},U_{\beta}]=2\pi i \delta_{\alpha \beta}U_{\beta}$$ 
This relation means that these operators specify a connection on the
module ${\cal S}$. It is equivalent to the Leibniz rule 
 $\nabla _{\alpha}(fe)=f\cdot \nabla _{\alpha}e+\delta _{\alpha}f\cdot e$, 
where $f=\sum c_{\vec {n}}U_{\vec {n}}\in T_{\theta}^d$,\ \ \  $e\in {\cal
S}$, and $\delta _{\alpha}$ stands for an infinitesimal automorphism
corresponding to translation in the direction $\alpha $. (An abelian   Lie
algebra $L={\bf R}^d$ acts naturally on the  algebra $T_{\theta}^d$ by
means of infinitesimal automorphisms - derivations). Generators of this 
Lie
algebra $\delta _1,...,\delta _d$ act in the following way: $\delta
_{\alpha} U_{\alpha}=2\pi i U_{\alpha}$ and $\delta _{\alpha}
U_{\beta} =0$ for $\alpha \not= \beta$.)

The curvature of the connection $\nabla _{\alpha}$ is
constant: $$F_{\alpha\beta}=[\nabla_{\alpha},\nabla _{\beta}]=2\pi
i\omega_{\alpha\beta}\cdot  1.$$

We say that a  noncommutative torus is equipped with complex structure  if
the Lie algebra $L={\bf R}^d$ acting on  $T_{\theta}^d$ is equipped with
such a structure. A complex structure on $L$ can be considered as a
decomposition of complexification $L\oplus iL$ of $L$ in  a direct sum of
two  complex  conjugate subspaces $L^{1,0}$ and $L^{0,1}$. We denote by
$\bar {\delta}_1,...,\bar {\delta}_d$ a basis in $L^{0,1}$; one can
express $\bar {\delta}_{\alpha}$ in terms of $\delta _{\alpha}$  as $\bar
{\delta}_{\alpha}=h_{\alpha}^{\beta}\delta _{\beta}$ where
$h_{\alpha}^{\beta}$ is a nondegenerate complex  $d\times d$ matrix. A
complex structure on a $T_{\theta}^d$-module  ${\cal E}$ can be 
defined as a collection of  ${\bf C}$-linear operators $\bar
{\nabla}_1,...,\bar {\nabla}_d$ on ${\cal E}$ satisfying 
\begin {equation}
\bar {\nabla}_{\alpha}(f\cdot e)=f\bar {\nabla}_{\alpha}e+\bar
{\delta}_{\alpha}f\cdot e,
\end {equation}

\begin {equation}
[\bar {\nabla}_{\alpha}, \bar {\nabla}_{\beta}]=0
\end {equation}
Here $f\in T_{\theta}^d,\ \ e\in {\cal E}$, Eqn (12) means that 
operators $\bar {\nabla}_1,...,\bar {\nabla}_d$ specify a 
 $\bar {\partial}$-connection  and Eqn (13) means that this $\bar
{\partial}$-connection is flat. 
A vector $e\in {\cal E}$ is called holomorphic if  $\bar
{\nabla}_{\alpha}e=0$, $\alpha =1,...,d$.

If   $ {\nabla}_{\alpha}$ is a constant curvature connection and $\bar
{\delta}_{\alpha}=h_{\alpha}^{\beta}\delta_{\beta}$ then 
\begin {equation}  
\bar {\nabla}_{\alpha}=h_{\alpha}^{\beta} {\nabla}_{\beta} 
\end {equation}
is a $\bar {\partial}$-connection; this connection is flat if $L^{0,1}$ is
an isotropic submanifold of $L\oplus iL$ with respect to (pre)symplectic
form specified by the curvature of connection 
${\nabla}_{\alpha}$. In particular, $T_{\theta}^d$-module  ${\cal S}$
constructed above by means of the matrix $2\pi \omega_{\alpha\beta}=2\pi
(\theta^{-1})_{\alpha\beta}$ can be equipped with a complex structure
corresponding to a Lagrangian submanifold $\Omega\subset L\oplus iL$. (We
consider $L\oplus iL$ as a symplectic manifold with symplectic structure
corresponding to the tensor $\omega _{\alpha\beta}.$)  The
$T_{\theta}^d$-module  ${\cal S}$ equipped with this complex structure
will be denoted by ${\cal S}_{\Omega}$. We will prove the following 
statement.

{\it  If} $\Omega$ {\it is a positive  Lagrangian submanifold of}
$L\oplus iL$ {\it then there exists a holomorphic vector in} ${\cal 
S}_{\Omega}$
{\it and this vector is unique up to a factor. We denote this vector by}
$\vartheta _{\Omega, \theta}$. {\it If} $\Omega$ {\it is not  positive 
then 
there
exist no non-zero  holomorphic vectors in} ${\cal S}_{\Omega}$.

Recall, that the tensor $\omega_{\alpha\beta}$ determines not only
bilinear symplectic pairing $<v,w>$ on $L\oplus iL$ that we always
consider, but also sesquilinear pairing $(v,w)=<\bar {v}, w>$.
One says that $\Omega$ is positive if the Hermitian form  $(v,v)$ has
positive imaginary part on $\Omega$ (i. e. ${\rm Im}(v,v)={\rm Im}<\bar
{v},v>$ is positive for $v\in \Omega,\ \  v\not= 0$).

To prove the theorem we work in coordinates $x^1,...,x^g$ where  the
symplectic form on $L$ is standard:  
$$\omega =dx^{\alpha}dp_{\alpha}$$
where $p_{\alpha}=x^{\alpha +g}$ and $\alpha =1,...,g$. Without loss of
generality we can assume that the natural projection $(x,p) \rightarrow p$
is surjective on Lagrangian subspace $\Omega$ (if this is not the case we
can change variables performing a symplectomorphism $\tilde
{x}^{\alpha}=p_{\alpha},\ \  \tilde {p}_{\alpha}=-x^{\beta}$ for several
indices $\alpha \in \{ 1,...,g\}$. This means that the equation of
$\Omega$ can be written in the form 
$$p_{\alpha}+\Omega_{\alpha\beta}x^{\beta}=0,\ \ \  \alpha=1,...,g$$ where
$\Omega_{\alpha\beta}$ is a symmetric complex matrix. The condition of
positivity of the Lagrangian subspace $\Omega$ means that the imaginary
part of the matrix $\Omega _{\alpha\beta}$ is positive. 

Representing the vector $\vartheta _{\Omega, \theta}$ as a function $\tau
(x^1,...,x^g)$ of variables $x^1,...,x^g$ we can write condition of
holomorphicity in the following form:
\begin {equation}
({\partial\over \partial
x^{\alpha}}-2\pi i\Omega_{\alpha\beta}x^{\beta})\tau=0.
\end {equation}
We see that up to a constant factor this function is a quadratic exponent:
\begin {equation}
\tau (x^1,...,x^n)=\exp(\pi ix^{\alpha}\Omega_{\alpha\beta} x^{\beta})
\end {equation}
In the case when $\Omega$ is a positive Lagrangian subspace we obtain that
$\tau$ belongs to the ${\cal S}({\bf R}^g)$ and therefore represents a
holomorphic vector in ${\cal S}_{\Omega}$  (theta-vector).  If $\Omega$ is
not positive the quadratic exponent (16) does not tend to zero at
infinity, hence does not belong to ${\cal S}({\bf R}^g)$; this means that
there are no holomorphic vectors in ${\cal S}_{\Omega}$. 

Let us present the above consideration in more invariant way. We can start
with $d$-dimensiomal vector space $L$ equipped with symplectic form
$\omega$. The invariant form of (18) is the relation
$$[\nabla_X, \nabla_Y]=2\pi i\omega_{XY}$$
where $\nabla _X$ is an anti-Hermitian operator in ${\cal S}$ that 
linearly
depends on $X\in L$. A theta-vector $\vartheta _{\Omega }$ in ${\cal S}$
is specified by means of positive Lagrangian $\Omega =L^{0,1}\subset
L\oplus iL$. Let us denote by $\gamma$ a linear operator on $L$ preserving
$\omega$. Then $\gamma \Omega =\tilde {\Omega}$ is again a positive
Lagrangian subspace of $L$; in other words the symplectic group Sp$(L)$
acts on Siegel upper half-space ${\cal H}_g$. It is well known that for
every element $\gamma\in {\rm Sp}(L)$ we can construct a unitary operator
$V_{\gamma}$ acting on ${\cal S}$ by means of the formula 
$$ V_{\gamma }\nabla _x V_{\gamma}^{-1}=\nabla _{\gamma x}$$
The operator $V_{\gamma}$ is defined up to a constant factor. The
correspondence $\gamma\rightarrow V_{\gamma}$ can be considered as a
projective representation of Sp$(L)$. It is evident that 
\begin {equation}
\vartheta _{\gamma \Omega }=V_{\gamma} \vartheta _{\Omega }. 
\end {equation}
(Recall that all objects entering (17) are defined up to a constant
factor.) 

Notice, that the definition of ${\cal S}$ was based only on symplectic
structure on $L$,  the definition of $\vartheta _{\Omega}$ on symplectic
and complex (=Kaehler) structure on $L$. To define a  structure of
$T_{\theta}^d$-module on ${\cal S}$ we should fix a basis $e_{\alpha}$ in
$L$; then the multiplicative generators of representation of
noncommutative torus $T_{\theta}^d$ with $\theta
^{\alpha\beta}=\omega(e_{\alpha},   e_{\beta})$ can be represented by
$U_{\alpha} =\exp (-\nabla_{e_{\alpha}})$.  If $\gamma\in {\rm  Sp}(L)$ 
then the operators $\tilde
{U}_{\alpha}=V_{\gamma}U_{\alpha}V_{\gamma}^{-1}$ specify an equivalent
representation of the same algebra $T_{\theta }^d$   (an isomorphic
$T_{\theta}^d$-module ). In other words, we can say that the group Sp$(L)$
acts on the space of complex structures on  $T_{\theta}^d$-module ${\cal
S}$ and that the transformation of theta-vectors $\vartheta
_{\theta,\Omega}$ by the action of $\gamma\in {\rm Sp}(L)$ is governed by
the formula (17).

 In the case when the entries of the matrix $\theta$ are integer numbers
the algebra $T_{\theta}^d$ is commutative; it is isomorphic to the algebra
of smooth functions on a torus and projective modules correspond to the
vector bundles. We will consider in detail the case when $\theta$ is the
standard symplectic matrix. We will show that in this case vectors
$\vartheta _{\Omega, \theta}$ are related to classical
theta-functions. Similar results can be proven for arbitraty matrix
$\theta$ with integer or, more generally, rational entries. 

In the case at hand we realize a projective module over $T_{\theta}^d$ as
the space ${\cal S}={\cal S}({\bf R}^g),\ \  d=2g$, with operators
$U_1,...,U_g$ acting as shifts:
$$(U_{\alpha}f)(\vec {x})=f(\vec {x}-\vec {e}_{\alpha})$$
and operators $U_{g+1},...,U_{2g}$ acting by the formula 
$$(U_{\alpha}f)(\vec {x})=e^{2\pi i x^{\alpha}}f(\vec {x}).$$
(Here $e_{\alpha}$ stands for the unit vector in the direction
$\alpha$.) This follows from (11) where covariant derivatives are taken in
the form (10).
$$(\nabla _{\alpha}f)(\vec {x})=  2\pi i \hat {x}^{\alpha}\ \ \  {\rm
for}\ 1\leq \alpha \leq g,$$
$$(\nabla _{\alpha}f)(\vec {x})= \partial _{\alpha}\ \ \  {\rm for}\ g<
\alpha \leq {2g}.$$
Unitary operators $U_1,...,U_{2g}$ commute, therefore there exists a
system $\varphi _{\vec {\rho},\vec {\sigma}}$ of common
(generalized) eigenfunctions of these operators. Namely, we can consider
distributions 
\begin {equation}
 \varphi _{\vec {\rho},\vec {\sigma}}(\vec {x})=\sum _{\vec {k}\in {\bf
Z}^g}e^{2\pi i \vec {\rho}\vec {k}}\delta (\vec {x}-\vec {\sigma}-\vec
{k}).
\end {equation}
It is easy to check that 
$$U_{\alpha}\varphi _{\vec {\rho},\vec {\sigma}}= e^{-2\pi
i\rho^{\alpha}}\varphi_{\vec {\rho} ,\vec {\sigma}}\ \ \  {\rm for}\ 1\leq
\alpha \leq g,$$
\begin {equation}
U_{\alpha}\varphi _{\vec {\rho},\vec {\sigma}}=e^{2\pi i\sigma
^{\alpha}}\varphi_{\vec {\rho} ,\vec {\sigma}}\ \ \  {\rm for}\ g< \alpha
\leq 2g,
\end {equation}
\begin {equation}
\varphi _{\vec {\rho}+\vec {e}, \vec {\sigma}}= \varphi _{\vec {\rho},\vec
{\sigma}},\ \ \      
\varphi _{\vec {\rho},\vec {\sigma}+\vec {e}}= e^{-2\pi i\vec {\rho}\vec
{e}}\varphi_{\vec {\rho} ,\vec {\sigma}}\ \ \  {\rm for\ \  every}\   \vec
{e}\in {\bf Z}^g.
\end {equation}
Using the functions $\varphi _{\vec {\rho},\vec {\sigma}}$ we can assign
to every vector $f\in {\cal S}$ a function
\begin {equation}
\tilde {f}(\vec {\rho},\vec {\sigma})=\int \bar {\varphi}_{\vec
{\rho},\vec {\sigma}}(\vec {x})f(\vec {x})dx=\sum  _{\vec {k}\in {\bf
Z}^g}e^{-2\pi  i\vec {\rho}\vec {k}}f (\vec {\sigma}+\vec {k}).
\end {equation}
This function obeys 
\begin {equation}
\tilde {f}(\vec {\rho}+\vec {e},\vec {\sigma})= \tilde {f}(\vec {\rho},
\vec {\sigma}), \ \ \  \tilde {f}(\vec {\rho}, \vec {\sigma}+\vec
{e})=e^{2\pi i\vec {\rho}\vec {e}} f(\vec {\rho} ,\vec {\sigma}) 
\end {equation}
and therefore can be considered as a section of a line bundle over a
torus; we will denote this bundle by $\Theta$. Conversely, having a
function $\tilde {f}$  obeying (22) we can restore $f$ using the formula 
\begin {equation}
f(\vec {x})=\int \tilde {f}_{\vec {\rho},\vec {\sigma}}\cdot \varphi
_{\vec {\rho},\vec {\sigma}}(\vec {x}) d \vec {\rho} d\vec {\sigma}
\end {equation}
Here we integrate over a fundamental domain of the lattice ${\bf Z}^d={\bf
Z}^g\times {\bf Z}^g$ acting on ${\bf R}^d={\bf R}^g\times {\bf
R}^g$. (The integrand is ${\bf Z}^d$-periodic, therefore the integral does
not depend on the choice of  fundamental domain.)

We obtained a correspondence between elements of projective module ${\cal
S}$ and sections of  
a line bundle $\Theta$ over  a torus. Covariant derivatives $\nabla
_{\alpha}$ correspond to the operators $\tilde {\nabla}_{\alpha}$ on
section of $\Theta $ that are given by the
formula 
$$\tilde {\nabla}_{\alpha}\tilde {f}_{\vec {\rho}, \vec {\sigma}}=(2\pi
i\sigma  ^{\alpha}-{\partial \over \partial \sigma^{\alpha}})\tilde
{f}_{\vec {\rho},\vec {\sigma}}\ \ \ {\rm for}\ \  1\leq \alpha\leq g,$$
\begin {equation}
\tilde {\nabla}_{\alpha}\tilde {f}_{\vec {\rho}, \vec {\sigma}}={\partial
\over \partial \sigma^{\alpha}}\tilde {f}_{\vec {\rho},\vec {\sigma}}\ \ \
{\rm for}\ \  g<\alpha\leq 2g.
\end {equation}
Let us suppose that a complex structure on $L$  and ${\cal S}$ is
constructed by means of Lagrangian subspace $L^{0,1}=\Omega \subset
L\oplus iL$ corresponding to the matrix $\Omega _{\alpha \beta}$.  Then
the functions 
\begin {equation}
z^i=(\Omega \vec {\sigma})^i -\rho ^i
\end {equation}
are holomorphic; they can be considered as complex coordinates on
$L$. Using these coordinates we can present the torus $T^d=L/{\bf Z}^d$ in
the form ${\bf C}^g/D$ where the lattice $D$ is spanned by vectors $\vec
{m}$ and vectors $\Omega \vec {m}$ with $\vec {m}\in {\bf Z}^g$. The
Lagrangian subspace $\Omega$ specifies a complex structure on the line
bundle $\Theta$; corresponding holomorphic section can be identified with
holomorphic function of $z^1,...,z^g$ obeying conditions (2),(3). More
precisely, if $f(x^1,...,x^n)$ represents a holomorphic vector in ${\cal
S}_{\Omega}$ then 
\begin {equation}
\tilde {f}_{\vec {\rho}, \vec {\sigma}}=e^{\pi i\vec {\sigma} \Omega
\vec {\sigma}}\Phi(\vec {z}),
\end {equation}
where $\Phi (z^1,...,z^g)$ is a holomorphic function obeying (2),(3).
We see that the vector $\vartheta _{\Omega, \theta}$ corresponds to the
theta-function (1) up to a factor.  This fact can be derived also from
(21), (25); such a derivation does not leave any unknown constant factors
in the relation between $\vartheta _{\Omega,\theta }$ and
theta-function. The calculation is trivial: it follows from (21) that the
function of $\vec {\rho}$ and $\vec {\sigma}$ that corresponds to the
vector (16) is equal to 
\begin {equation}
\tilde {\tau}_{\vec {\rho},\vec {\sigma}}=\sum  _{\vec {k}\in {\bf
Z}^g}e^{-2\pi  i\vec {\rho}\vec {k}}\exp  (\pi i(\vec {\sigma}+\vec
{k})\Omega (\vec {\sigma}+\vec {k}))
\end {equation}
Expressing $\tilde {\tau}$ in terms of $z^1,...,z^g$ we obtain
\begin {equation}
\tilde {\tau}_{\vec {\rho},\vec {\sigma}}=e^{\pi i
\vec {\sigma}\Omega\vec {\sigma}}\vartheta _{\Omega}(\vec {z})
\end {equation}
in agreement with (25).

Using the above statements it is easy to analyze modular properties of
theta-functions with respect to ${\rm Sp}(2g,{\bf Z})$
transformations. The behavior  of theta-vector  with respect to ${\rm
Sp}(L)$ is given by the formula (17). To analyze the transformation  of
theta-function with respect to $\gamma \in {\rm Sp}(L)$ we should know
also the behavior of $\varphi_{\vec {\rho}, \vec {\sigma}}$; corresponding
formulas can be obtained easily if $\gamma \in {\rm Sp}(2g,{\bf
Z})$. Recall that $\varphi_{\vec {\rho}, \vec {\sigma}}$ were defined as
(generalized) eigenvectors of $U_{\alpha}=\exp (-\theta
^{\alpha\beta}\nabla_{\beta})$. Transformed functions $\varphi_{\vec
{\rho}, \vec {\sigma}}$  can be considered as eigenvectors of
$U_{\alpha}^{\prime}=\exp (-\theta ^{\alpha\beta}\nabla_{\beta}^{\prime})$
where 
$$\nabla_{\alpha}^{\prime}= V_{\gamma} \nabla_{\alpha} 
V_{\gamma}^{-1}=\gamma_{\alpha}^{\beta}
\nabla _{\beta}$$ 
and  $\gamma = (\gamma _{\beta}^{\alpha})\in {\rm Sp}(L)$ is a symplectic
matrix with integer entries. It is easy to see that the operators
$U_{\alpha}^{\prime}$ belong to the algebra generated by commuting
operators $U_{\alpha}$ and therefore  $\varphi_{\vec {\rho}, \vec
{\sigma}}$ is an eigenvector of $U_{\alpha}^{\prime}$. More precisely, 
$$U_{\alpha}^{\prime}=\exp (-\theta
^{\alpha\beta}\nabla_{\beta}^{\prime})=\exp (-\theta ^{\alpha\beta}\gamma
_{\beta}^{\lambda}\nabla_{\lambda})=\exp (-\gamma _{\lambda}
^{\alpha}\theta ^{\lambda\mu}\nabla_{\mu})=c_{\alpha}\Pi_{1\leq \lambda
\leq 2g}U_{\lambda}^{\gamma_{\lambda}^{\alpha}},$$  where 
$$c _{\alpha}=\exp (\pi i \gamma _{\lambda}^{\alpha}\theta
^{\lambda\mu}\gamma_{\mu}^{\alpha}).$$ 
This formula follows immediately from the relation 
$$ \exp (K_1+...+K_n)=\exp ({1\over 2}\sum _{i<j}[K_i,K_j])\cdot
\Pi_{1\leq i\leq n}e^{K_i}$$
that is valid in the case when all commutators $[K_i,K_j]$ commute with
$K_1,...,K_n$. It is easy to check that $c_{\alpha}=1$ iff $\gamma\in
\Gamma _{1,2}$. Using this fact we can obtain the well known modular
properties of theta-functions. ( A similar proof of these properties is 
given in [10].)

There exist interesting generalizations of the notion of theta-function
and corresponding generalizations of the notion of theta-vector. In these
generalizations theta-functions are defined as some holomorphic sections
of vector bundles over tori and theta-vectors as holomorphic vectors in
projective modules over noncommutative tori. We will consider basic
modules over noncommutative tori [7]. These modules can be characterized
as projective $T_{\theta}$-modules  with constant curvature connections
that have another  noncommutative torus $T_{\hat {\theta}}$ as an
endomorphism algebra:  $T_{\hat {\theta}}={\rm End}_{T_{\theta}}{\cal
E}$. (Basic modules can be described also in terms of their $K$-theory
classes; see [7].) Basic $T_{\theta}$-module ${\cal E}$ can be considered
as $(T_{\theta},T_{\hat {\theta}})$-bimodule that  establishes
(gauge) Morita equivalence [6] between $T_{\theta}$ and $T_{\hat
{\theta}}$; antisymmetric matrices $\theta$ and $\hat {\theta}$ are
related by means of 
fractional linear transformation $\hat 
{\theta}=(M\theta
+N)(R\theta +S)^{-1}$ governed by an element of $SO(d,d,{\bf Z})$. 
We will
denote a constant curvature connection on basic $T_{\theta}$-modules
${\cal E}$ by $\nabla_{\alpha}$ and assume that its 
curvature  $[\nabla_{\alpha},\nabla _{\beta}]=2\pi i \omega _{\alpha
\beta}$ is a non-degenerate matrix. 

Then, as we have seen, a complex structure on ${\cal E}$ corresponds to a
Lagrangian subspace $\Omega$ of $L\oplus iL$. We define the space ${\cal
H}_{\Omega}$ of theta-vectors as the subspace of ${\cal E}$ consisting of
holomorphic vectors with respect to complex structure on ${\cal E}$
corresponding to $\Omega$.  The space ${\cal H}_{\Omega}$ can be easily
described by means of results of [8] where we gave a construction of
${\cal E}$ and analyzed the moduli space of 
constant curvature connections on ${\cal E}$. Namely, ${\cal E}$ can be
regarded as a direct sum of $N$ copies of irreducible representation of
canonical commutation relations (8).  {\it If  Lagrangian subspace} 
$\Omega$  {\it is
positive then the dimension of} ${\cal H}_{\Omega}$ {\it is equal to} $N$; 
{\it if}
$\Omega$ {\it is not positive theta-vectors don't exist.} 

The number $N$ can be expressed in terms of topological numbers (of
$K$-theory class)  of the module ${\cal E}$. Namely, the $K$-theory
class $\mu ({\cal E})$ of ${\cal E}$ can be regarded as an integer element
of Grassmann algebra with $2g$ generators $\alpha ^1,...,\alpha^{2g}$ (or
as a sequence of integer antisymmetric tensors).The $K$-theory class
of a basic module is a generalized quadratic exponent ( a limit of
quadratic exponents) having relatively prime coefficients. {\it It
follows from the
formula (18) of [8] that} 
$$\mu ({\cal E})=N\alpha ^1\alpha ^2...\alpha^{2g}+{\rm lower\ \  order \
\  terms}$$
{\it (i.e.} $N$ {\it is determined by the antisymmetric tensor of highest 
possible
rank).} Using Berezin integration on Grassmann algebra we can say 
that $N$ is an integral of $\mu({\cal E})$ with respect to anticommuting 
variables $\alpha ^1,...,\alpha ^{2g}$.

 Notice that using the connection 
$\nabla _{\alpha}$ we can identify
the Lie algebra $L$ acting on $T_{\theta}$ and  the Lie algebra $\hat {L}$
acting on $T_{\hat {\theta}}$. This means that {\it a complex structure 
on}
$T_{\theta}$-{\it module} ${\cal E}$  {\it induces a complex structure on} 
${\cal E}$
{\it considered as} 
$T_{\hat {\theta}}$-{\it module} (in other words we can regard ${\cal E}$ 
as
complex  
$(T_{\theta}, T_{\hat {\theta}})$-bimodule.  The notion of theta-vector in
complex $T_{\theta}$-module ${\cal E}$ coincides with the notion of
theta-vector in corresponding $ T_{\hat {\theta}}$-module. If  we have
along with a basic $(T_{\theta}, T_{\hat {\theta}})$-bimodule  ${\cal E}$
 a basic  $(T_{\hat {\theta}}, T_{\hat {\hat {\theta}}})$-bimodule
${\cal E}^{\prime}$  we can consider the tensor product 
$${\cal E}^{\prime \prime}={\cal E}\otimes {\cal E}^{\prime}$$
as a basic  $(T_{\theta}, T_{\hat {\hat {\theta}}})$-bimodule.   If tori
$T_{\theta}, T_{\hat {\theta}}$       and  $T_{\hat {\hat {\theta}}}$as
well as
bimodules ${\cal E}, {\cal E}^{\prime}$ are equipped with  complex
structures then it is easy to construct a  complex structure on $(T_
{\theta}, T_{\hat {\hat {\theta}}})$-bimodule ${\cal E}^{\prime \prime}$. 

One can check that {\it this construction determines a map of the tensor
product} ${\cal H}\otimes {\cal H}^{\prime}$ {\it into} ${\cal H}^{\prime
\prime}$. Here ${\cal H},{\cal H}^{\prime}, {\cal H}^{\prime \prime}$
stand for spaces of holomorphic vectors (theta-vectors) in ${\cal E},
{\cal E}^{\prime}$ and ${\cal E}^{\prime \prime}$ correspondingly. 

The tensor product ${\cal E}\otimes {\cal E}^{\prime}$ can be described 
explicitly by means of results of [6], [8]. It follows from these results 
that a triple $(\theta, \hat {\theta},g)$ where $\theta, \hat {\theta}$ 
are antisymmetric $d\times d$ matrices,  $g\in SO(d,d,{\bf Z})$ and $\hat 
{\theta}=\tilde {g}\theta$ determines a basic $(T_{\theta},T_{\hat 
{\theta}})$-bimodule ${\cal E}$;  all basic bimodules correspond to such 
triples.  (Here $\tilde {g}$ stands for fractional linear 
transformation corresponding to 
 $g\in SO(d,d,{\bf Z})$.)

 {\it The tensor product } 
${\cal E}^{\prime \prime}$ {\it of bimodules} ${\cal E}$ {\it and} ${\cal 
E}^{\prime}$ {\it corresponding to triples} 
$(\theta, \hat {\theta},g)$ 
{\it and} $(\theta^{\prime}, \hat {\theta}^{\prime},g^{\prime})${\it is 
well defined in the case when } $\hat {\theta}=\theta ^{\prime}$; {\it it 
corresponds to a triple} $(\theta, \hat {\theta}^{\prime},g^{\prime}g)$. 
An explicit description of the space of theta-vectors and of a map ${\cal 
H}\otimes {\cal H}^{\prime}\rightarrow {\cal H}^{\prime\prime}$ can be 
obtained by means of results of [6], [8]. In particular, the dimension of 
the space ${\cal H}$ of theta-vectors in basic $(T_{\theta}, T_{\hat 
{\theta}})$-bimodule ${\cal E}$, corresponding a triple 
$(\theta,\hat{\theta},g),$ can be expressed in terms of Berezin integral 
of $\mu=S(g)\cdot 1$. (The group $SO(d,d,{\bf Z})$ acts on Grassmann 
algebra generated by $\alpha^1,...,\alpha ^d$ by means of spinor 
representation $S(g)$.) To verify this statement we consider ${\cal E}$ as 
Morita equivalence bimodule; then ${\cal E}$ regarded as $T_{\hat 
{\theta}}$-module, corresponds to $T_{\theta}^1$ (to free 
$T_{\theta}$-module of rank $1$).

It seems that the theory of theta-vectors should be closely related to
Manin's theory of quantized theta-functions [11]-[14], in particular, the
partial multiplication constructed by  Manin, should correspond to
the tensor product of bimodules.  It would be interesting to analyze this
relation thoroughly.

  We have seen that noncommutative analogs of theta-functions are related 
to projective $T_{\theta}$-modules equipped with complex structure. 
Complex projective modules appear also in the study on noncommutative 
instantons (see [15]). We expect that in the framework of Connes' 
noncommutative geometry one can develop complex noncommutative geometry, 
that is related in some way to noncommutative algebraic geometry. This is 
an interesting open problem. 
 
 \centerline {\bf Acknowledgements.}

I am indebted to Yu. Manin and Ya. Soibelman for interesting discussions.

\centerline {\bf References.}

1. A. Connes, {\it  C*-algebres et Geometrie Differentielle}, C. R.
Acad. Sci. Paris {\bf 290} (1980), 599-604.           
  
  2. A. Connes and  M. Rieffel, {\it Yang-Mills for  Noncommutative
Two-Tori}, Contemporary Math. {\bf 66} (1987), 237-266.

   3. M. Rieffel, {\it Projective Modules over Higher-dimensional
Noncommutative Tori}, Can. J. Math., Vol. {\bf XL}, No. 2 (1988), 257-338.     

 4. M. Rieffel  and  A. Schwarz, {\it Morita Equivalence of
Multidimensional Noncommutative Tori},  Intl. J. of Math  {\bf 10(2)}
(1999), 289-299                   
 
5. A. Schwarz, {\it Gauge theories on noncommutative spaces},  ICMP
lecture,
hep-th/0011261. 
 
6. A.  Schwarz, {\it  Morita Equivalence and Duality},  Nucl. Phys.
{\bf  B 534} (1998), 720-738.
                                                     
  7  . A. Konechny  and  A.  Schwarz, {\it  BPS States on Noncommutative
Tori and Duality}, Nucl. Phys. {\bf B 550}  (1999), 561-584.       

  8. A. Konechny  and  A. Schwarz,  {\it  Moduli Spaces of Maximally
Supersymmetric Solutions on Noncommutative Tori and  Noncommutative
Orbifolds},  JHEP {\bf 09} (2000), 1-23. 

  9. D. Mumford, {\it Tata Lectures on Theta I}, Birkhauser (1983).

  10. D. Mumford, {\it Tata Lectures on Theta III} (with M. 
Nori and P. Norman), Birkhauser (1991).
 
  11. Yu. Manin, {\it Quantized theta-functions}, In Common Trends in
Mathematics and Quantum Field Theories (Kioto,1990), Progress of
Theor. Phys. Supplement {\bf 102 } (1990), 219-228.           

  12. Yu. Manin, {\it Mirror symmetry and quantization of abelian
varieties},  preprint math. AG/0005143.

  13. Yu. Manin, {\it Theta Functions, Quantum Tori and Heisenberg
Groups}, preprint math.AG/0011197 

  14. Ya. Soibelman, {\it Quantum tori, mirror symmetry and deformation 
theory}, preprint math.QA/0011162 

  15. A.  Schwarz, {\it Noncommutative instantons: a new approach}, 
hep-th/0102182
 
\end{document}